\documentclass[letterpaper, 12pt]{amsart}

\usepackage[all]{xy}
\usepackage{latexsym}
\usepackage{amssymb}
\usepackage{fullpage}
\usepackage{amsmath}
\usepackage{amsthm}
\usepackage{amscd}
\usepackage{epic,eepic}
\usepackage[dvips]{graphicx}
\usepackage{yfonts}

\newtheorem{theorem}{Theorem}[section]
\newtheorem{lemma}[theorem]{Lemma}
\newtheorem{corollary}[theorem]{Corollary}
\newtheorem{definition}[theorem]{Definition}
\newtheorem{proposition}[theorem]{Proposition}

\numberwithin{equation}{section}

\newcommand{\V}{\Vert}

\newcommand{\la}{\langle}
\newcommand{\ra}{\rangle}
\newcommand{\ga}{\mathfrak{A}}
\title{A Noncommutative Gauss Map}
\author{Caleb Eckhardt}
\address{Department of Mathematics, University of Illinois, Urbana, IL, 61801}
\email[Caleb Eckhardt]{ceckhard@gmail.edu}
\date{}

\begin{document}
\maketitle
\begin{abstract} The aim of this paper is to transfer the Gauss map, which is a Bernoulli shift for continued fractions, to the noncommutative setting.  We feel that a natural place for such a map to act is on the AF algebra $\ga$  considered separately by F. Boca and D. Mundici.  The center of $\ga$ is isomorphic to $C[0,1]$, so we first consider the action of the Gauss map on $C[0,1]$ and then extend the map to $\ga$ and show that the extension inherits many desirable properties. 
\end{abstract}

\section{Introduction and Notation}

Florin Boca in \cite{Boca08} and Daniele Mundici in \cite{Mundici88} separately considered an AF algebra $\ga$ that is associated with the Farey tessellation. The algebra $\ga$ exhibits many interesting properties, not the least of which is the connection between $\ga$ and the unit interval $[0,1].$  This connection is not merely topological, but also number theoretic. We briefly explain this connection 

Let $Z(\ga)$ denote the center of $\ga.$  As noted in \cite{Boca08}, we have $C[0,1]\cong Z(\ga).$  Moreover, the maximal ideal space of $\ga$ is homeomorphic (when equipped with the topology induced by Prim($\ga$)) to $[0,1]$ in a natural way \cite[Corollary 12]{Boca08}.  For each irrational $0<\theta<1$, let $\mathcal{J}_\theta$ denote the maximal ideal of $\ga$ associated to $\theta.$  It was shown in \cite{Boca08} that $\ga/\mathcal{J}_\theta\cong \mathfrak{F}_\theta$, the Effros-Shen algebra, defined in \cite{Effros80}, associated with the continued fraction expansion of $\theta.$  

In other words, if we employ the topological decomposition theory of $C^*$-algebras and visualize $\ga$ as continuous, operator-valued functions on its maximal ideal space (a visualization which is usually ``incorrect, but fruitful'' \cite[Page 91]{Pedersen79}), then each function evaluated at $\theta$ takes values in the Effros-Shen algebra $\mathfrak{F}_\theta.$ So it is not simply the topology around $\theta$ that determines this visualization, but also the continued fraction expansion of $\theta.$

Given the close connection between $\ga$ and the continued fraction expansions of numbers in $[0,1]$, it is natural to try and extend important functions from number theory (especially those related to continued fractions) to the $C^*$-algebra $\ga.$ The Gauss map might be the most fundamental such function, hence we take it as our starting point  

Recall the Gauss map $G:[0,1]\rightarrow [0,1]$  defined by 
$G(0)=0$ and $G(x)=1/x-\lfloor 1/x \rfloor$ if $x\neq0,$ where $\lfloor\cdot\rfloor $ denotes the greatest integer function.
One can think of $G$ as the Bernoulli shift for continued fractions.  Indeed, given $\theta\in[0,1]$ with continued fraction expansion $\theta=[a_1,a_2,...]$ then $G([a_1,a_2,...])=[a_2,a_3,...].$ One can also recover the continued fraction expansion of $\theta$ by implementations of $G$ and $\lfloor\cdot\rfloor.$

 We first consider the induced action of $G$ on $Z(\ga)=C[0,1].$
 First note that $f\circ G\in C[0,1]$ if and only if $f$ is a constant function.  
Therefore we will consider the ``adjoint'' action of the Gauss map on $C[0,1].$ Let $\mu$ denote Gauss measure on $[0,1]$ defined by $d\mu=\frac{d\theta}{\ln 2 (\theta+1)}$, where $d\theta$ denotes Lebesgue measure.  Then $G$ is $\mu$-invariant, i.e. $\mu(G^{-1}(E))=\mu(E)$ for every Borel set $E\subseteq [0,1]$ (see \cite{Iosifescu02} for details). From this it follows that the map
\begin{equation*}
 V_G(f)(\theta)=f(G(\theta))\quad \textrm{ for }f\in L^2(\mu), \theta\in[0,1]
\end{equation*}
is an isometry.   A standard calculation reveals that
\begin{equation}
 V_G^*(f)(\theta)=\sum_{s=1}^\infty f\Big( \frac{1}{\theta+s}\Big) \frac{1+\theta}{(\theta+s)(\theta+s+1)} \label{eq:V_G*definition}
\end{equation}
and it is routine to verify that $V_G^*(f)\in C[0,1]$  when $f\in C[0,1].$  We mention that, symbolically, $V_G^*$ is the Perron-Frobenius operator of $G$ under $\mu$ and refer the reader to \cite[Chapter 2]{Iosifescu02} for details about Perron-Frobenius operators and their connections to continued fractions.

Furthermore, if we embed $C[0,1]$ into $B(L^2(\mu))$ as $f\mapsto \mathcal{M}_f$ where $\mathcal{M}_f(g)=fg$ for all $g\in L^2(\mu)$, then
\begin{equation}
 V_G^* \mathcal{M}_f V_G=\mathcal{M}_{V_G^*(f)}.  \label{eq:conjugationproperty}
\end{equation}
This defines a unital completely positive map, which we will henceforth denote by $\mathbb{G}$, on $C[0,1]$ that not only respects the action of $G$ on its maximal ideal space, i.e. for each $E\subseteq [0,1]$  let $J_E$ denote the ideal of $C[0,1]$ consisting of those functions that vanish on $E$, then
\begin{equation}
 \mathbb{G}(J_{G^{-1}(E)})\subseteq J_{E}, \label{eq:mappingGaussideals}
\end{equation}
but is also $\mu-$invariant, i.e.
\begin{equation}
 \int f d\mu=\int \mathbb{G}(f)d\mu\quad\textrm{ for all }\quad f\in C[0,1]. \label{eq:muinvariant}
\end{equation}
Hence we are looking for an extension of $\mathbb{G}$ to $\ga$ that satisfies the natural analogs of (\ref{eq:conjugationproperty})-(\ref{eq:muinvariant}).  In order to do this we must first consider what Gauss measure should mean on $\ga.$ Our first step is proving that every state on $C[0,1]$ has a unique extension to a trace on $\ga$ (Theorem \ref{thm:statetrace}).
 D. Mundici showed \cite[Theorem 4.5]{Mundici08} that the state space of $C[0,1]$ and the space of tracial states on $\ga$ are affinely, weak* homeomorphic.  But for our purposes, we will need the extension property from Theorem \ref{thm:statetrace}. 

For reasons that will become clear, we have to slightly modify the natural analogs of (\ref{eq:conjugationproperty}) and (\ref{eq:muinvariant}).
In particular, we will  use Theorem \ref{thm:statetrace} to define two separate state extensions, $\phi$ and $\tau$, of $\mu$ and intertwine between these two GNS representations to obtain analogs of (\ref{eq:conjugationproperty}) and (\ref{eq:muinvariant}).  Let $(\pi_\phi, L^2(\ga,\phi))$ and $(\pi_\tau, L^2(\ga,\tau))$  be the GNS representations of $\ga$ associated with $\phi$ and $\tau.$  Since $\phi$ and $\tau$ are extensions of $\mu$, it follows that $L^2(\mu)\subseteq L^2(\ga,\phi), L^2(\ga,\tau)$ and 
\begin{equation*}
 \pi_\phi(f)|_{L^2(\mu)}=\mathcal{M}_f \quad \textrm{ for every }f\in Z(\ga)\cong C[0,1].
\end{equation*}
This allows us to prove the main theorem:

\begin{theorem} \label{thm:maintheorem} There is  a unital completely positive map $\widetilde{\mathbb{G}}:\ga\rightarrow\ga$ and an isometry $\widetilde{V}_G: L^2(\ga,\tau)\rightarrow L^2(\ga,\phi)$ such that
 \begin{enumerate}
 \item $\widetilde{\mathbb{G}}|_{C[0,1]}=\mathbb{G}.$ 
 \item $\widetilde{\mathbb{G}}(\mathcal{J}(G^{-1}(E)))\subseteq \mathcal{J}(E),$ for each  $E\subset [0,1]$ \textup{(}$\mathcal{J}_E$ are the ideals of $\ga$ defined in \cite{Boca08}\textup{)}.
  
  \item $\widetilde{V}_G|_{L^2(\mu)}=V_G$ and $\widetilde{V}_G^*|_{L^2(\mu)}=V_G^*.$
  \item $\widetilde{V}_G^*\pi_\phi(x)\widetilde{V}_G=\pi_\tau(\widetilde{\mathbb{G}}(x))$ for $x\in\ga.$ Hence $\widetilde{V}_G^*\pi_\phi(f)\widetilde{V}_G|_{L^2(\mu)}=\mathcal{M}_{\mathbb{G}(f)}$ for $f\in C[0,1].$
  \item $\phi(x)=\tau(\widetilde{\mathbb{G}}(x))$  for $x\in\ga.$
 \end{enumerate}
\end{theorem}

In order to set our notation, we now recall some relevant facts about the AF algebra $\mathfrak{A}$ defined in \cite{Boca08} and \cite{Mundici88}. We will use the same notation as in \cite{Boca08}, in particular $p(n,k), q(n,k)\in\mathbb{Z}^+$ and $r(n,k)=\frac{p(n,k)}{q(n,k)}$ for $n\geq0$ and $0\leq k\leq 2^n$ all have the same meaning  and we will frequently refer to the relationships between them as defined on \cite[pg. 3]{Boca08}. Recall that $\ga$ is the inductive limit of the finite dimensional $C^*$-algebras,
\begin{equation*}
 \mathfrak{A}_n=\bigoplus_{0\leq k\leq 2^n} M_{q(n,k)}.
\end{equation*}
For the convenience of the reader, and with thanks to F. Boca for supplying us with the code, we reproduce the Bratelli diagram of $\ga$ from \cite[Figure 2]{Boca08}.
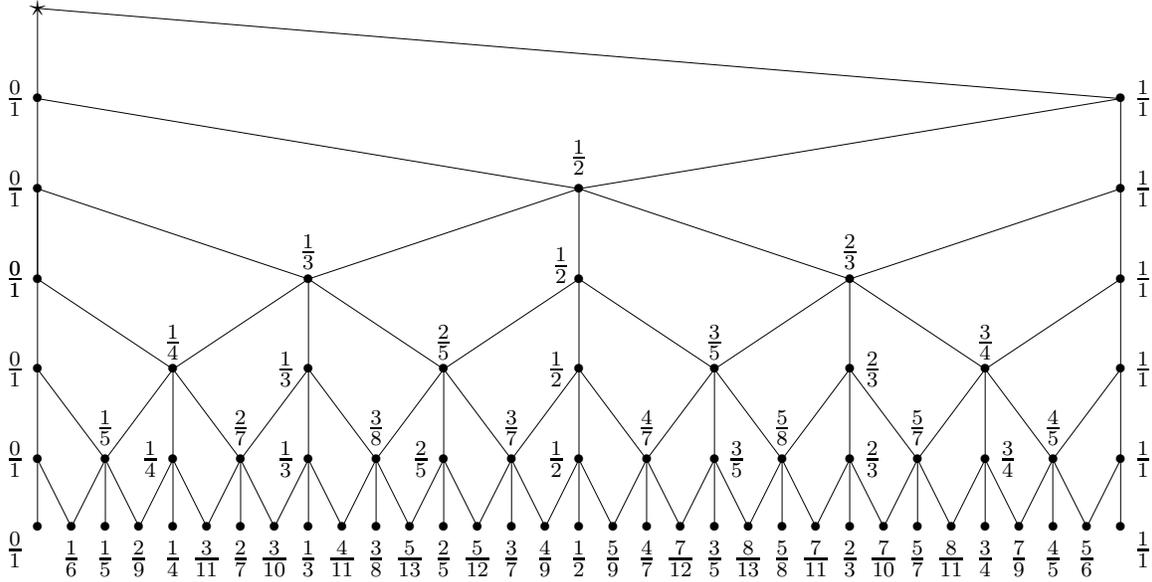
\begin{figure}[htb]
\begin{center}
\unitlength 0.6mm
\begin{picture}(0,120)(0,0)
%\thinlines
\path(-120,100)(-120,120)(120,100)
\path(-120,80)(-120,100)(0,80)(120,100)(120,80)
\path(-120,80)(-120,60)
\path(-120,60)(-120,80)(-60,60)(0,80)(60,60)(120,80)(120,60)
\path(0,80)(0,60)

\path(-120,40)(-120,60)(-90,40)(-60,60)(-30,40)(0,60)(30,40)(60,60)(90,40)(120,60)(120,40)
\path(-60,60)(-60,40) \path(0,60)(0,40) \path(60,60)(60,40)
\path(-120,20)(-120,40)(-105,20)(-90,40)(-75,20)(-60,40)(-45,20)(-30,40)(-15,20)(0,40)(15,20)
(30,40)(45,20)(60,40)(75,20)(90,40)(105,20)(120,40)(120,20)
\path(-90,40)(-90,20) \path(-60,40)(-60,20) \path(-30,40)(-30,20)
\path(0,40)(0,20) \path(30,40)(30,20) \path(60,40)(60,20)
\path(90,40)(90,20)

\path(-120,5)(-120,20)(-112.5,5)(-105,20)(-97.5,5)(-90,20)(-82.5,5)(-75,20)(-67.5,5)(-60,20)
(-52.5,5)(-45,20)(-37.5,5)(-30,20)(-22.5,5)(-15,20)(-7.5,5)(0,20)
\path(120,5)(120,20)(112.5,5)(105,20)(97.5,5)(90,20)(82.5,5)(75,20)(67.5,5)(60,20)
(52.5,5)(45,20)(37.5,5)(30,20)(22.5,5)(15,20)(7.5,5)(0,20)

\path(-105,20)(-105,5) \path(-90,20)(-90,5) \path(-75,20)(-75,5)
\path(-60,20)(-60,5) \path(-45,20)(-45,5) \path(-30,20)(-30,5)
\path(-15,20)(-15,5) \path(0,20)(0,5) \path(105,20)(105,5)
\path(90,20)(90,5) \path(75,20)(75,5) \path(60,20)(60,5)
\path(45,20)(45,5) \path(30,20)(30,5) \path(15,20)(15,5)

\put(-125,100){\makebox(0,0){{\small $\frac{0}{1}$}}}
\put(0,87){\makebox(0,0){{\small $\frac{1}{2}$}}}
\put(125,100){\makebox(0,0){{\small $\frac{1}{1}$}}}
\put(-125,80){\makebox(0,0){{\small $\frac{0}{1}$}}}
\put(125,80){\makebox(0,0){{\small $\frac{1}{1}$}}}
\put(-125,60){\makebox(0,0){{\small $\frac{0}{1}$}}}
\put(-60,66){\makebox(0,0){{\small $\frac{1}{3}$}}}
\put(-4,63){\makebox(0,0){{\small $\frac{1}{2}$}}}
\put(60,66){\makebox(0,0){{\small $\frac{2}{3}$}}}
\put(125,60){\makebox(0,0){{\small $\frac{1}{1}$}}}
\put(-125,60){\makebox(0,0){{\small $\frac{0}{1}$}}}

\put(-125,40){\makebox(0,0){{\small $\frac{0}{1}$}}}
\put(-90,46){\makebox(0,0){{\small $\frac{1}{4}$}}}
\put(-65,40){\makebox(0,0){{\small $\frac{1}{3}$}}}
\put(-30,46){\makebox(0,0){{\small $\frac{2}{5}$}}}
\put(-5,40){\makebox(0,0){{\small $\frac{1}{2}$}}}
\put(30,46){\makebox(0,0){{\small $\frac{3}{5}$}}}
\put(65,40){\makebox(0,0){{\small $\frac{2}{3}$}}}
\put(90,46){\makebox(0,0){{\small $\frac{3}{4}$}}}
\put(125,40){\makebox(0,0){{\small $\frac{1}{1}$}}}

\put(-120,100){\makebox(0,0){{\tiny $\bullet$}}}
\put(120,100){\makebox(0,0){{\tiny $\bullet$}}}
\put(-120,80){\makebox(0,0){{\tiny $\bullet$}}}
\put(0,80){\makebox(0,0){{\tiny $\bullet$}}}
\put(120,80){\makebox(0,0){{\tiny $\bullet$}}}
\put(-120,60){\makebox(0,0){{\tiny $\bullet$}}}
\put(-60,60){\makebox(0,0){{\tiny $\bullet$}}}
\put(0,60){\makebox(0,0){{\tiny $\bullet$}}}
\put(60,60){\makebox(0,0){{\tiny $\bullet$}}}
\put(120,60){\makebox(0,0){{\tiny $\bullet$}}}

\put(-120,40){\makebox(0,0){{\tiny $\bullet$}}}
\put(-90,40){\makebox(0,0){{\tiny $\bullet$}}}
\put(-60,40){\makebox(0,0){{\tiny $\bullet$}}}
\put(-30,40){\makebox(0,0){{\tiny $\bullet$}}}
\put(0,40){\makebox(0,0){{\tiny $\bullet$}}}
\put(30,40){\makebox(0,0){{\tiny $\bullet$}}}
\put(60,40){\makebox(0,0){{\tiny $\bullet$}}}
\put(90,40){\makebox(0,0){{\tiny $\bullet$}}}
\put(120,40){\makebox(0,0){{\tiny $\bullet$}}}

\put(-125,20){\makebox(0,0){{\small $\frac{0}{1}$}}}
\put(-105,27){\makebox(0,0){{\small $\frac{1}{5}$}}}
\put(-95,20){\makebox(0,0){{\small $\frac{1}{4}$}}}
\put(-75,27){\makebox(0,0){{\small $\frac{2}{7}$}}}
\put(-65,20){\makebox(0,0){{\small $\frac{1}{3}$}}}
\put(-45,27){\makebox(0,0){{\small $\frac{3}{8}$}}}
\put(-35,20){\makebox(0,0){{\small $\frac{2}{5}$}}}
\put(-15,27){\makebox(0,0){{\small $\frac{3}{7}$}}}
\put(-5,20){\makebox(0,0){{\small $\frac{1}{2}$}}}
\put(125,20){\makebox(0,0){{\small $\frac{1}{1}$}}}
\put(105,27){\makebox(0,0){{\small $\frac{4}{5}$}}}
\put(95,20){\makebox(0,0){{\small $\frac{3}{4}$}}}
\put(75,27){\makebox(0,0){{\small $\frac{5}{7}$}}}
\put(65,20){\makebox(0,0){{\small $\frac{2}{3}$}}}
\put(45,27){\makebox(0,0){{\small $\frac{5}{8}$}}}
\put(35,20){\makebox(0,0){{\small $\frac{3}{5}$}}}
\put(15,27){\makebox(0,0){{\small $\frac{4}{7}$}}}
\put(-120,120){\makebox(0,0){$\star$}}

\put(-120,20){\makebox(0,0){{\tiny $\bullet$}}}
\put(-105,20){\makebox(0,0){{\tiny $\bullet$}}}
\put(-90,20){\makebox(0,0){{\tiny $\bullet$}}}
\put(-75,20){\makebox(0,0){{\tiny $\bullet$}}}
\put(-60,20){\makebox(0,0){{\tiny $\bullet$}}}
\put(-45,20){\makebox(0,0){{\tiny $\bullet$}}}
\put(-30,20){\makebox(0,0){{\tiny $\bullet$}}}
\put(-15,20){\makebox(0,0){{\tiny $\bullet$}}}
\put(0,20){\makebox(0,0){{\tiny $\bullet$}}}
\put(120,20){\makebox(0,0){{\tiny $\bullet$}}}
\put(105,20){\makebox(0,0){{\tiny $\bullet$}}}
\put(90,20){\makebox(0,0){{\tiny $\bullet$}}}
\put(75,20){\makebox(0,0){{\tiny $\bullet$}}}
\put(60,20){\makebox(0,0){{\tiny $\bullet$}}}
\put(45,20){\makebox(0,0){{\tiny $\bullet$}}}
\put(30,20){\makebox(0,0){{\tiny $\bullet$}}}
\put(15,20){\makebox(0,0){{\tiny $\bullet$}}}

\put(-125,0){\makebox(0,0){{\small $\frac{0}{1}$}}}
\put(-112.5,-2){\makebox(0,0){{\small $\frac{1}{6}$}}}
\put(-105,-2){\makebox(0,0){{\small $\frac{1}{5}$}}}
\put(-97.5,-2){\makebox(0,0){{\small $\frac{2}{9}$}}}
\put(-90,-2){\makebox(0,0){{\small $\frac{1}{4}$}}}
\put(-82.5,-2){\makebox(0,0){{\small $\frac{3}{11}$}}}
\put(-75,-2){\makebox(0,0){{\small $\frac{2}{7}$}}}
\put(-67.5,-2){\makebox(0,0){{\small $\frac{3}{10}$}}}
\put(-60,-2){\makebox(0,0){{\small $\frac{1}{3}$}}}
\put(-52.5,-2){\makebox(0,0){{\small $\frac{4}{11}$}}}
\put(-45,-2){\makebox(0,0){{\small $\frac{3}{8}$}}}
\put(-37.5,-2){\makebox(0,0){{\small $\frac{5}{13}$}}}
\put(-30,-2){\makebox(0,0){{\small $\frac{2}{5}$}}}
\put(-22.5,-2){\makebox(0,0){{\small $\frac{5}{12}$}}}
\put(-15,-2){\makebox(0,0){{\small $\frac{3}{7}$}}}
\put(-7.5,-2){\makebox(0,0){{\small $\frac{4}{9}$}}}
\put(0,-2){\makebox(0,0){{\small $\frac{1}{2}$}}}
\put(125,0){\makebox(0,0){{\small $\frac{1}{1}$}}}
\put(112.5,-2){\makebox(0,0){{\small $\frac{5}{6}$}}}
\put(105,-2){\makebox(0,0){{\small $\frac{4}{5}$}}}
\put(97.5,-2){\makebox(0,0){{\small $\frac{7}{9}$}}}
\put(90,-2){\makebox(0,0){{\small $\frac{3}{4}$}}}
\put(82.5,-2){\makebox(0,0){{\small $\frac{8}{11}$}}}
\put(75,-2){\makebox(0,0){{\small $\frac{5}{7}$}}}
\put(67.5,-2){\makebox(0,0){{\small $\frac{7}{10}$}}}
\put(60,-2){\makebox(0,0){{\small $\frac{2}{3}$}}}
\put(52.5,-2){\makebox(0,0){{\small $\frac{7}{11}$}}}
\put(45,-2){\makebox(0,0){{\small $\frac{5}{8}$}}}
\put(37.5,-2){\makebox(0,0){{\small $\frac{8}{13}$}}}
\put(30,-2){\makebox(0,0){{\small $\frac{3}{5}$}}}
\put(22.5,-2){\makebox(0,0){{\small $\frac{7}{12}$}}}
\put(15,-2){\makebox(0,0){{\small $\frac{4}{7}$}}}
\put(7.5,-2){\makebox(0,0){{\small $\frac{5}{9}$}}}

\put(-120,5){\makebox(0,0){{\tiny $\bullet$}}}
\put(-112.5,5){\makebox(0,0){{\tiny $\bullet$}}}
\put(-105,5){\makebox(0,0){{\tiny $\bullet$}}}
\put(-97.5,5){\makebox(0,0){{\tiny $\bullet$}}}
\put(-90,5){\makebox(0,0){{\tiny $\bullet$}}}
\put(-82.5,5){\makebox(0,0){{\tiny $\bullet$}}}
\put(-75,5){\makebox(0,0){{\tiny $\bullet$}}}
\put(-67.5,5){\makebox(0,0){{\tiny $\bullet$}}}
\put(-60,5){\makebox(0,0){{\tiny $\bullet$}}}
\put(-52.5,5){\makebox(0,0){{\tiny $\bullet$}}}
\put(-45,5){\makebox(0,0){{\tiny $\bullet$}}}
\put(-37.5,5){\makebox(0,0){{\tiny $\bullet$}}}
\put(-30,5){\makebox(0,0){{\tiny $\bullet$}}}
\put(-22.5,5){\makebox(0,0){{\tiny $\bullet$}}}
\put(-15,5){\makebox(0,0){{\tiny $\bullet$}}}
\put(-7.5,5){\makebox(0,0){{\tiny $\bullet$}}}
\put(0,5){\makebox(0,0){{\tiny $\bullet$}}}
\put(120,5){\makebox(0,0){{\tiny $\bullet$}}}
\put(112.5,5){\makebox(0,0){{\tiny $\bullet$}}}
\put(105,5){\makebox(0,0){{\tiny $\bullet$}}}
\put(97.5,5){\makebox(0,0){{\tiny $\bullet$}}}
\put(90,5){\makebox(0,0){{\tiny $\bullet$}}}
\put(82.5,5){\makebox(0,0){{\tiny $\bullet$}}}
\put(75,5){\makebox(0,0){{\tiny $\bullet$}}}
\put(67.5,5){\makebox(0,0){{\tiny $\bullet$}}}
\put(60,5){\makebox(0,0){{\tiny $\bullet$}}}
\put(52.5,5){\makebox(0,0){{\tiny $\bullet$}}}
\put(45,5){\makebox(0,0){{\tiny $\bullet$}}}
\put(37.5,5){\makebox(0,0){{\tiny $\bullet$}}}
\put(30,5){\makebox(0,0){{\tiny $\bullet$}}}
\put(22.5,5){\makebox(0,0){{\tiny $\bullet$}}}
\put(15,5){\makebox(0,0){{\tiny $\bullet$}}}
\put(7.5,5){\makebox(0,0){{\tiny $\bullet$}}}

\end{picture}
\end{center}
\caption{Bratelli diagram of $\ga$} \label{Figure1}
\end{figure}

For each $n\geq0$  let $\mathbb{E}_n:\ga\rightarrow \ga_n$ be conditional expectations such that 
\begin{equation}
 \mathbb{E}_n\mathbb{E}_m=\mathbb{E}_m\mathbb{E}_n \quad \textrm{ for all }n,m\geq0. \label{eq:commutingexpectations}
\end{equation}
The existence of such conditional expectations is guaranteed by Arveson's extension theorem, or since $\ga$ is AF, one can construct such maps explicitly.  Furthermore for $0\leq k\leq 2^n$ let  $\mathbb{E}_{(n,k)}:\ga\rightarrow M_{q(n,k)}\subseteq \ga_n$ be conditional expectations, such that 
\begin{equation}
\mathbb{E}_{(n,k)}\mathbb{E}_n=\mathbb{E}_n\mathbb{E}_{(n,k)}\quad \textrm{for }n\geq0\quad \textrm{and }0\leq k\leq 2^n. \label{eq:commutingfdexpectations} 
\end{equation}
Note that for each $x\in\ga$ we have
\begin{equation}
 \lim_{n\rightarrow\infty}\mathbb{E}_n(x)=x. \label{eq:convergencethroughexpectations}
\end{equation}
We will use the following notation throughout:
 For a unital $C^*$-algebra $A$, we let 
\begin{enumerate}
\item[$\cdot$] $Z(A)$ denote the center of $A$, 
\item[$\cdot$] $\mathcal{S}(A)$ denote the state space of $A$,
\item[$\cdot$] $\mathcal{T}(A)$ denote the set of all unital traces of $A.$
\item[$\cdot$] $M_n$ denote $n\times n$ matrices over $\mathbb{C}$,
\item[$\cdot$] $\tau_n$ the unital trace on $M_n$ and $1_n\in M_n$ the identity.
\item[$\cdot$] For $\frac{p}{q}\in\mathbb{Q}\cap[0,1]$ in reduced form we define
\begin{equation*}
 M_{\frac{p}{q}}:=M_q\quad  \tau_{\frac{p}{q}}:=\tau_q \in \mathcal{T}(M_q).
\end{equation*}
\end{enumerate}

\section{State Extensions and Conditional Expectation onto $Z(\ga)$} \label{sec:center}
In this section we will construct a conditional expectation from $\ga$ onto $Z(\ga)$ that preserves every $\tau\in\mathcal{T}(\ga).$ This will provide the key step in the proof of Theorem \ref{thm:statetrace}.

\begin{definition} \label{defn:rationaltraces} Let $n\geq0$ and $0\leq k\leq 2^n.$  Define $\tau_{(n,k)}\in\mathcal{T}(\ga)$ as
\begin{equation*}
 \tau_{(n,k)}(x)=\tau_{q(n,k)}\circ \mathbb{E}_{(n,k)}(x).
\end{equation*}
\end{definition}
The following lemma is immediate from (\ref{eq:commutingexpectations}) and (\ref{eq:commutingfdexpectations}):
\begin{lemma} \label{lem:welldefinedf_x}
 For $n\geq0$, $0\leq k\leq 2^n$ and   $\ell\geq0$, we have
\begin{equation*}
 \tau_{(n,k)}=\tau_{(n+\ell,2^\ell k)}.
\end{equation*}
\end{lemma}

\begin{proposition} Let $x\in\ga.$ Define the function $f_x:\mathbb{Q}\cap[0,1]\rightarrow\mathbb{C}$ as
\begin{equation}
 f_x(r(n,k))=\tau_{(n,k)}(x)\quad \textrm{for }n\geq0\textrm{ and } 0\leq k\leq 2^n.  \label{eq:defnoff_x}
\end{equation}
Then $f_x$ is well-defined and extends to a continuous function on $[0,1].$  
\end{proposition}
\begin{proof} If $r(n,k)=r(n',k')$ with $n'\geq n$ then there is an $\ell\geq0$ such that $n'=n+\ell$ and $k'=2^\ell k.$
Hence $f_x$ is well-defined by Lemma \ref{lem:welldefinedf_x}.

Let $n\geq0$ and $0<k<2^n.$ By the relationships defined in \cite[pg. 3]{Boca08}, the following function is continuous and  piecewise affine on $[0,1]$: 
 \begin{equation*}
 B_{(n,k)}(\theta)=\left\{ \begin{array}{ll}0 & \textrm{ if }0\leq \theta\leq r(n,k-1)\\
q(n,k)\Big(q(n,k-1)\theta-p(n,k-1)\Big) & \textrm{ if }r(n,k-1)\leq \theta\leq r(n,k)\\
q(n,k)\Big(p(n,k+1)-q(n,k+1)\theta\Big) &\textrm{ if }r(n,k)\leq \theta\leq r(n,k+1)\\
0& \textrm{ if }r(n,k+1)\leq \theta\leq 1
\end{array} \right. 
\end{equation*}

We first let $x\in\ga_n\subset\ga$ and prove that $f_x$ extends to a continuous function on $[0,1].$ 
Suppose first that $0<2k+1<2^n$ and $\mathbb{E}_{(n,2k+1)}(x)=x.$  Without loss of generality suppose that $\tau_{(n,2k+1)}(x)=1.$  We show that $f_x=B_{(n,2k+1)}.$

It is clear that $f_x|_{[0,r(n,2k)]\cup[r(n,2k+2),1]}\equiv0.$  We now show by induction on $\ell\geq0$ that
\begin{equation}
 (\forall \ell\geq0)(\forall\textrm{ } 2^\ell 2k< j < 2^\ell(2k+2))(f_x(r(n+\ell,j))=B_{(n,2k+1)}(r(n+\ell,j))).\label{eq:inductionhypo}
\end{equation}
For $\ell=0$, we have $f_x(r(n,2k+1))=\tau_{(n,2k+1)}(x)=1=B_{(n,2k+1)}(r(n,2k+1)).$  Suppose now that (\ref{eq:inductionhypo}) holds for $\ell\geq0$ and prove (\ref{eq:inductionhypo}) for $\ell+1.$

If $j=2i$ is even,  then
\begin{align*}
 &f_x(r(n+\ell+1,2i))=\tau_{(n+\ell+1,2i)}(x)=\tau_{(n+\ell,i)}(x)=f_x(r(n+\ell,i))\\
=&B_{(n,2k+1)}(r(n+\ell,i))=B_{(n,2k+1)}(r(n+\ell+1,2i)).
\end{align*}
If $j=2i+1$ is odd, then
\begin{align*}
 &f_x(r(n+\ell+1,2i+1))=\tau_{(n+\ell+1,2i+1)}(x)\\
&=\frac{q(n+\ell,i)}{q(n+\ell+1,2i+1)}\tau_{(n+\ell,i)}(x)+\frac{q(n+\ell,i+1)}{q(n+\ell+1,2i+1)}\tau_{(n+\ell,i+1)}(x)\\
&=\frac{q(n+\ell,i)}{q(n+\ell+1,2i+1)}B_{(n,2k+1)}(r(n+\ell,i))+\frac{q(n+\ell,i+1)}{q(n+\ell+1,2i+1)}B_{(n,2k+1)}(r(n+\ell,i+1))\\
&=B_{(n,2k+1)}(r(n+\ell+1,2i+1)).
\end{align*}
Here the last line follows by the relationships in \cite[pg. 3]{Boca08} and because $B_{(n,2k+1)}$ is piecewise affine.
This shows that (\ref{eq:inductionhypo}) holds, hence $f_x$ extends to a continuous function on $[0,1].$

Now suppose that $0<2^mk<2^n$ with $k$ odd and $\mathbb{E}_{(n,2^mk)}(x)=x.$  Then, 
\begin{equation*}
 x=\mathbb{E}_{(n-m,k)}(x)-\mathbb{E}_{(n-m+1,2k-1)}(x)-\mathbb{E}_{(n-m+1,2k+1)}(x).
\end{equation*}
So, by the first part of the proof it follows that $f_x$ is continuous.

For  $x=1\oplus 0\oplus \cdots\oplus 0, 0\oplus \cdots \oplus 0\oplus 1\in \ga_n$, the proof that $f_x$ is continuous is exactly the same as above, so we omit the proof.  This shows that for every $n\geq0$ and each $x\in\ga_n$ that $f_x$ is continuous.  Moreover note that the linear map $x\mapsto f_x$ defined on $\bigcup_{n=1}^\infty \ga_n$ is contractive, hence $f_x$ is continuous for every $x\in\ga.$
\end{proof}

As observed in \cite{Boca08}, $Z(\ga)\cong C[0,1].$ We now construct an explicit isomorphism.
For each $n\geq0,$ define $\mathcal{Z}_n:C[0,1]\rightarrow Z(\ga_n)\subset \ga$ by
\begin{equation}
\mathcal{Z}_n(f)=\bigoplus_{0\leq k\leq 2^n} f(r(n,k))1_{q(n,k)}  \label{eq:defofZ_n}
\end{equation}
By \cite[pg. 3]{Boca08}, for each $n\geq0$ we have $\max\{|r(n,k)-r(n,k+1)|:0\leq k<2^n\}=1/(n+1).$  Hence for $m\geq n$ we have
\begin{equation*}
 \V \mathcal{Z}_n(f)-\mathcal{Z}_m(f)\V\leq \sup\{ |f(\theta)-f(\theta')|:|\theta-\theta'|\leq 1/(n+1)\}.
\end{equation*}
Therefore $\mathcal{Z}_n(f)$ is a Cauchy sequence in $\ga$ because $f$ is uniformly continuous on $[0,1].$

Define $\mathcal{Z}:C[0,1]\rightarrow \ga$ by
\begin{equation}
 \mathcal{Z}(f)=\lim_{n\rightarrow\infty} \mathcal{Z}_n(f).\label{eq:Tdefn}
\end{equation}

\begin{theorem} \label{thm:conditionalexpectationdefn} The map $\mathcal{Z}:C[0,1]\rightarrow Z(\ga)$ is a *-isomorphism.  Moreover the map
\newline $\mathbb{E}_Z:\ga\rightarrow Z(\ga)$ defined by
\begin{equation*}
 \mathbb{E}_Z(x)=\mathcal{Z}(f_x)
\end{equation*}
is a conditional expectation such that 
\begin{equation}
\tau(\mathbb{E}_Z(x))=\tau(x)\quad \textrm{for every }\quad \tau\in\mathcal{T}(\ga).\label{eq:tracepreserving}
\end{equation}
\end{theorem}
\begin{proof}
By (\ref{eq:Tdefn}) it is clear that $\mathcal{Z}$ is  a *-monomorphism, and since  $\mathcal{Z}_n(f)\in Z(\ga_n)$ for each $n\geq0$, it follows that $\mathcal{Z}(f)\in Z(\ga).$  We now show that $\mathcal{Z}$ is surjective.
Let $n\geq0$ and $y\in \ga_n.$  Then 
\begin{equation}
y\in Z(\ga_n)\quad \textrm{if and only if}\quad y=\bigoplus_{0\leq k\leq 2^n}\tau_{(n,k)}(y)1_{q(n,k)}.\label{eq:charofcenterinfdalgebras}
\end{equation}
Let $x\in Z(\ga).$  By (\ref{eq:defnoff_x}) and (\ref{eq:defofZ_n}) it follows that
\begin{equation*}
\mathcal{Z}_n(f_x)=\bigoplus_{0\leq k\leq 2^n}\tau_{(n,k)}(x)1_{q(n,k)}\in Z(\ga_n).
\end{equation*}
Since $x\in Z(\ga)$, it follows from (\ref{eq:convergencethroughexpectations}) that 
\begin{equation*}
 \lim_{n\rightarrow\infty} dist(\mathbb{E}_n(x), Z(\ga_n))=0,
\end{equation*}
from which we deduce by (\ref{eq:charofcenterinfdalgebras}) that $\mathcal{Z}_n(f_x)\rightarrow x.$
Therefore
\begin{equation}
 \mathbb{E}_Z(x)=\mathcal{Z}(f_x)=\lim_{n\rightarrow\infty}\mathcal{Z}_n(f_x)= x. \label{eq:surjectivityofZ}
\end{equation}
This shows that $\mathcal{Z}$ is surjective and also that $\mathbb{E}_Z$ is a conditional expectation.
We now show that $\mathbb{E}_Z$ preserves every trace of $\ga$. 
Let $\tau\in\mathcal{T}(\ga).$  By (\ref{eq:convergencethroughexpectations}) it follows that 
$\tau$ is the weak*-limit of $\tau\circ \mathbb{E}_n.$  Since $\tau\circ \mathbb{E}_n|_{\ga_n}\in\mathcal{T}(\ga_n)$, there is a convex combination of scalars $(\lambda_{(n,k)})_{0\leq k\leq 2^n}$ such that 
\begin{equation*}
 \tau\circ\mathbb{E}_n=\sum_{0\leq k\leq 2^n}\lambda_{(n,k)}\tau_{(n,k)}.
\end{equation*}
It follows that $\mathcal{T}(\ga)$ equals the weak* closure of the convex hull of the set $\{\tau_{(n,k)}:n\geq0, 0\leq k\leq 2^n\}.$  Therefore, we only need to check (\ref{eq:tracepreserving}) for the traces $\tau_{(n,k)}.$  To this end, let $x\in\ga$ then
\begin{equation}
 \tau_{(n,k)}(x)=f_x(r(n,k))=\tau_{(n,k)}(\mathcal{Z}_n(f_x))=\tau_{(n,k)}(\mathcal{Z}(f_x))=\tau_{(n,k)}(\mathbb{E}_Z(x)).  \label{eq:diracpointmassesabound}
\end{equation}
\end{proof}

\begin{theorem} \label{thm:statetrace}
The restriction map $\tau\mapsto \tau|_{Z(\ga)}$ defines a weak* homeomorphism from $\mathcal{T}(\ga)$ onto $\mathcal{S}(C[0,1]).$ 
In particular, every state on $Z(\mathfrak{A})$ has a unique tracial extension to $\mathfrak{A}.$
\end{theorem}
\begin{proof} Injectivity and weak*-continuity of the inverse both follow from (\ref{eq:tracepreserving}).  By (\ref{eq:diracpointmassesabound}) it follows that the restriction of $\tau_{(n,k)}$ to $Z(\ga)\cong C[0,1]$ is the Dirac measure $\delta_{\{r(n,k)\}},$  which shows surjectivity.
\end{proof}

\section{Ideals of $\ga$ and traces of $\ga$} \label{sec:ideals}
\begin{definition}
 Fix $\theta\in [0,1].$   Define $\tau^\ga_\theta\in \mathcal{T}(\ga)$ as the unique tracial extension of the Dirac measure $\delta_{\{\theta\}}\in \mathcal{S}(C[0,1])$ given by Theorem \ref{thm:statetrace}.  Note that for each $n\geq0$ and $0\leq k\leq 2^n$, we have $\tau_{r(n,k)}^\ga=\tau_{(n,k)}$ from Definition \ref{defn:rationaltraces}.
\end{definition}

For each $\theta\in[0,1]$, we recall the maximal ideals $I_\theta\subset \ga$ defined in \cite[Proposition 4]{Boca08}.
The following is a consequence of the proof of \cite[Proposition 4]{Boca08} and the correspondence made in Theorem \ref{thm:statetrace}.
\begin{corollary} \label{cor:traceideal} Fix $\theta\in [0,1].$   Then
\begin{equation}
 I_\theta=\{ x\in\ga: \tau_\theta^\ga(x^*x)=0\}. \label{eq:idealsandtraces}
\end{equation} 
\end{corollary}
Fix $\frac{p(n,k)}{q(n,k)}=\frac{p}{q}\in\mathbb{Q}\cap (0,1)$ in reduced form.  We define the *-homomorphism 
\begin{equation}
\pi_{\frac{p}{q}}:\ga\rightarrow M_{\frac{p}{q}} \label{eq:rationalhomo}
\end{equation}
 as ``evaluation along the path $r(n,k), r(n+1, 2k),..., r(n+\ell, 2^\ell k), ...$ in the Bratteli diagram.'' In particular, $\textrm{ker}(\pi_{\frac{p}{q}})=I_{\frac{p}{q}}$ (see  \cite[Proposition 4.(ii)]{Boca08} for details). We note that
\begin{equation}
 \tau^\ga_{\frac{p}{q}}(x)=\tau_{\frac{p}{q}}(\pi_{\frac{p}{q}}(x))\quad \textrm{ for every }x\in\ga.  \label{eq:tracesthroughrepresentations}
\end{equation}

\section{Construction of $\widetilde{\mathbb{G}}$}
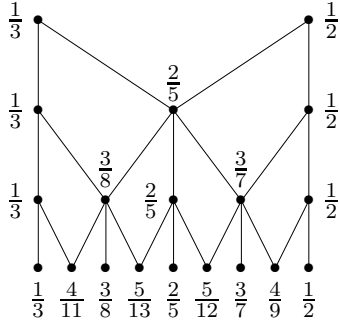
\begin{figure}[htb]
\begin{center}
\unitlength 0.6mm
\begin{picture}(-30,60)(0,0)
%\thinlines

\path(-60,60)(-30,40)(0,60)
\path(-60,60)(-60,40) \path(0,60)(0,40) 
\path(-60,40)(-45,20)(-30,40)(-15,20)(0,40)
 \path(-60,40)(-60,20) \path(-30,40)(-30,20)
\path(0,40)(0,20) 

\path(-60,20)
(-52.5,5)(-45,20)(-37.5,5)(-30,20)(-22.5,5)(-15,20)(-7.5,5)(0,20)

\path(-60,20)(-60,5) \path(-45,20)(-45,5) \path(-30,20)(-30,5)
\path(-15,20)(-15,5) \path(0,20)(0,5)

\put(-65,60){\makebox(0,0){{\small $\frac{1}{3}$}}}
\put(5,60){\makebox(0,0){{\small $\frac{1}{2}$}}}

\put(-65,40){\makebox(0,0){{\small $\frac{1}{3}$}}}
\put(-30,46){\makebox(0,0){{\small $\frac{2}{5}$}}}
\put(5,40){\makebox(0,0){{\small $\frac{1}{2}$}}}

\put(-60,60){\makebox(0,0){{\tiny $\bullet$}}}
\put(0,60){\makebox(0,0){{\tiny $\bullet$}}}

\put(-60,40){\makebox(0,0){{\tiny $\bullet$}}}
\put(-30,40){\makebox(0,0){{\tiny $\bullet$}}}
\put(0,40){\makebox(0,0){{\tiny $\bullet$}}}

\put(-65,20){\makebox(0,0){{\small $\frac{1}{3}$}}}
\put(-45,27){\makebox(0,0){{\small $\frac{3}{8}$}}}
\put(-35,20){\makebox(0,0){{\small $\frac{2}{5}$}}}
\put(-15,27){\makebox(0,0){{\small $\frac{3}{7}$}}}
\put(5,20){\makebox(0,0){{\small $\frac{1}{2}$}}}

\put(-60,20){\makebox(0,0){{\tiny $\bullet$}}}
\put(-45,20){\makebox(0,0){{\tiny $\bullet$}}}
\put(-30,20){\makebox(0,0){{\tiny $\bullet$}}}
\put(-15,20){\makebox(0,0){{\tiny $\bullet$}}}
\put(0,20){\makebox(0,0){{\tiny $\bullet$}}}

\put(-60,-2){\makebox(0,0){{\small $\frac{1}{3}$}}}
\put(-52.5,-2){\makebox(0,0){{\small $\frac{4}{11}$}}}
\put(-45,-2){\makebox(0,0){{\small $\frac{3}{8}$}}}
\put(-37.5,-2){\makebox(0,0){{\small $\frac{5}{13}$}}}
\put(-30,-2){\makebox(0,0){{\small $\frac{2}{5}$}}}
\put(-22.5,-2){\makebox(0,0){{\small $\frac{5}{12}$}}}
\put(-15,-2){\makebox(0,0){{\small $\frac{3}{7}$}}}
\put(-7.5,-2){\makebox(0,0){{\small $\frac{4}{9}$}}}
\put(0,-2){\makebox(0,0){{\small $\frac{1}{2}$}}}

\put(-60,5){\makebox(0,0){{\tiny $\bullet$}}}
\put(-52.5,5){\makebox(0,0){{\tiny $\bullet$}}}
\put(-45,5){\makebox(0,0){{\tiny $\bullet$}}}
\put(-37.5,5){\makebox(0,0){{\tiny $\bullet$}}}
\put(-30,5){\makebox(0,0){{\tiny $\bullet$}}}
\put(-22.5,5){\makebox(0,0){{\tiny $\bullet$}}}
\put(-15,5){\makebox(0,0){{\tiny $\bullet$}}}
\put(-7.5,5){\makebox(0,0){{\tiny $\bullet$}}}
\put(0,5){\makebox(0,0){{\tiny $\bullet$}}}

\end{picture}
\end{center}
\caption{Bratelli diagram of $\ga/\mathcal{J}_2$} \label{Figure2}
\end{figure}

In this section we construct our noncommutative Gauss map $\widetilde{\mathbb{G}}:\ga\rightarrow \ga.$
 Let $s\geq1.$  As in \cite[(3.1)]{Boca08} we define
\begin{equation*}
 \mathcal{J}_s:=\mathcal{J}([\frac{1}{s+1},\frac{1}{s}])=\bigcap_{\theta\in [\frac{1}{s+1},\frac{1}{s}]}I_\theta.
\end{equation*}
By Theorem \ref{thm:statetrace} and Section \ref{sec:ideals} we have
\begin{equation}
 \mathcal{J}_s=\textrm{ker}\Big(\bigoplus_{\frac{1}{s+1}<\frac{p}{q}\in\mathbb{Q}<\frac{1}{s}}\pi_{\frac{p}{q}}\Big). \label{eq:anotherdescriptionJ_s}
\end{equation}
For each $s\geq1$ the Bratelli diagram of $\ga/\mathcal{J}_s$ is the subdiagram of the Bratelli diagram of $\ga$ obtained by deleting all of the nodes 
\begin{equation*}
\{ r(n,k):r(n,k)\not\in [1/(s+1),1/s]\}\cup\{r(n,k):n<s\},
\end{equation*}
and deleting all edges connected to any of these nodes.  See Figure \ref{Figure2} for the Bratelli diagram of $\ga/\mathcal{J}_2.$

For each $s\geq1$, define the homeomorphism $g_s:[0,1]\rightarrow [1/(s+1),1/s]$ as 
\begin{equation*}
 g_s(\theta)=\frac{1}{\theta+s}, \label{eq:defng_s}
\end{equation*}
and recall that these maps are the building blocks for the commutative Gauss map $\mathbb{G}:C[0,1]\rightarrow C[0,1]$ defined in (\ref{eq:V_G*definition}).  Then consider the induced isomorphism 
\begin{equation*}
(g_s)_*:C[1/(s+1),1/s]\rightarrow C[0,1]\quad \textrm{defined by}\quad (g_s)_*(f)=f\circ g_s.
\end{equation*}
Since our goal is to extend $\mathbb{G}$ to a map on $\ga$, we first consider extensions of the maps $(g_s)_*$ as maps from $\ga/\mathcal{J}_s$ into $\ga.$
Unfortunately, there is no hope for these extensions to also be isomorphisms.  Indeed, by considering the Bratelli diagrams of $\ga$ and $\ga/\mathcal{J}_s$ it is clear that $K_0(\ga)\cong K_0(\ga/\mathcal{J}_s)$, but there is no unit-preserving, positive homomorphism that implements this isomorphism.  Hence $\ga\not\cong \ga/\mathcal{J}_s.$ We do the next best thing by defining a (non-unital) *-monomorphism $H_s:\ga\rightarrow \ga/\mathcal{J}_s$ and a unital completely positive (UCP for short) map $G_s:\ga/\mathcal{J}_s\rightarrow\ga$ such that $G_sH_s=id_\ga,$ and such that $G_s$ is an extension of $(g_s)_*.$ More importantly, the maps $G_s$ and $H_s$ will provide a nice relationship (see (\ref{eq:maintracerelationshipall})) between $\mathcal{T}(\ga)$ and $\mathcal{T}(\ga/\mathcal{J}_s).$  

For $n\geq0, $ let  $A_n\in  M_{2^{n+1}+1,2^n+1}(\mathbb{Z}^+)$  be the connecting homomorphisms from $\ga_n$ into $\ga_{n+1}$ such that
\begin{equation}
 \ga=\varinjlim(\ga_n, A_n) \label{eq:Fareyasinductivelimit}
\end{equation}
For example we have,
\begin{equation*}
 A_0=\left[ \begin{array}{cc} 1&0\\ 1&1\\ 0&1\\ \end{array}\right]\in M_{3,2}\textrm{ },\quad A_1=\left[ \begin{array}{ccc} 1&0&0\\ 1&1&0 \\ 0&1&0\\ 0&1&1\\ 0&0&1\\ \end{array}\right]\in M_{5,3}\textrm{ }, \cdots
\end{equation*}

For $n\geq0$, we define
\begin{equation*}
 (\ga/\mathcal{J}_s)_n:= \bigoplus_{0\leq k\leq 2^n}M_{\frac{q(n,k)}{p(n,k)+sq(n,k)}}=\bigoplus_{0\leq k\leq 2^n}M_{g_s(r(n,k))}.
\end{equation*}
By the description of the Bratteli diagram of $\ga/\mathcal{J}_s$ (see also Figure \ref{Figure2}) given above it follows that
\begin{equation}
 \ga/\mathcal{J}_s=\varinjlim((\ga/\mathcal{J}_s)_n, A_n) \label{eq:quotientasinductivelimit}
\end{equation}

Let $\ell^\infty(s)$ denote the $s$ dimensional, commutative $C^*$-algebra.  Consider the $C^*$-algebra,
\begin{equation*}
 \ell^\infty(s)\otimes \ga=\varinjlim (\ell^\infty(s)\otimes \ga_n, id_{\ell^\infty(s)}\otimes A_n)
\end{equation*}
Define $S=\left[ \begin{array}{llll} 1&1&\cdot&1 \end{array}\right]\in M_{1,s}.$ It is easy to see (using only the fact that $A_n\in M_{2^{n+1}+1,2^n+1}(\mathbb{Z}^+)$) that
\begin{equation*}
 A_n(S\otimes 1_{2^n+1})=(S\otimes 1_{2^{n+1}+1})1_s\otimes A_n\quad \textrm{ for every }n\geq0.
\end{equation*}

Hence, for each $n\geq0$ we are able to  define a *-homomorphism
\begin{equation*}
 \sigma_n: \ell^\infty(s)\otimes \ga_n=\bigoplus_{0\leq k\leq 2^n}\ell^\infty(s)\otimes M_{\frac{p(n,k)}{q(n,k)}}\rightarrow \bigoplus_{0\leq k\leq 2^n}M_{\frac{q(n,k)}{p(n,k)+sq(n,k)}}=(\ga/\mathcal{J}_s)_n
\end{equation*}
given by the matrix $S\otimes 1_{2^n+1}\in M_{2^n+1, s(2^n+1)}$ such that the following diagram commutes for every $n\geq0$
\begin{equation}
 \xymatrix{   (\ga/\mathcal{J}_s)_n \ar[r]^{A_n} & (\ga/\mathcal{J}_s)_{n+1} \\
\ell^\infty(s)\otimes \ga_n \ar[u]^{\sigma_n} \ar[r]^{1_s\otimes A_n} & \ell^\infty(s)\otimes \ga_{n+1} \ar[u]^{\sigma_{n+1}}  }  \label{eq:commutingsquare1}
\end{equation}
Let $e_1,...,e_s\in\ell^\infty(s)$ denote the standard basis. 
For each $n\geq0$ define the UCP map $V_n:(\ga/\mathcal{J})_n\rightarrow \ell^\infty(s)\otimes \ga_n$ by
\begin{equation*}
 V_n(x)=\sigma_n^{-1}\Big( \sum_{i=1}^s \sigma_n(e_i\otimes 1_{\ga_n})x\sigma_n(e_i\otimes 1_{\ga_n})\Big)
\end{equation*}
Define $\psi_s\in \mathcal{S}(\ell^\infty(s))$ by
\begin{equation*}
 \psi_s\Big(\sum_{i=1}^s \alpha_ie_i\Big)=\frac{1}{s}\sum_{i=1}^s \alpha_i.
\end{equation*}
It now follows from (\ref{eq:commutingsquare1})  that the following diagram commutes for all $n\geq0:$

\begin{equation}
 \xymatrix @R=.7in @C=2in    {   (\ga/\mathcal{J}_s)_n \ar[r]^{A_n} \ar@<1ex>[d]^{V_n} & (\ga/\mathcal{J}_s)_{n+1} \ar@<1ex>[d]^{V_{n+1}} \\
\ell^\infty(s)\otimes \ga_n \ar@<1ex>[u]^{\sigma_n} \ar[r]^{1_s\otimes A_n} \ar@<1ex>[d]^{\psi_s\otimes id_{\ga_n}} & \ell^\infty(s)\otimes \ga_{n+1} \ar@<1ex>[u]^{\sigma_{n+1}} \ar@<1ex>[d]^{\psi_s\otimes id_{\ga_{n+1}}}\\
\ga_n \ar@<1ex>[u]^{1_s\otimes id_{\ga_n}} \ar[r]^{A_n} & \ga_{n+1} \ar@<1ex>[u]^{1_s\otimes id_{\ga_{n+1}}} }   \label{eq:commutingsquare2}
\end{equation}
Furthermore,
\begin{equation}
 (\psi_s\otimes id_{\ga_n})\circ V_n\circ \sigma_n\circ(1_s\otimes id_{\ga_n})=id_{\ga_n}\quad \textrm{ for all }n\geq0. \label{eq:paritalinversesonfd}
\end{equation}

Now, let $x\in M_{\frac{p(n,k)}{q(n,k)}}\subset \ga_n$  and $y\in M_{\frac{q(n,k)}{p(n,k)+sq(n,k)}}\subset (\ga/\mathcal{J}_s)_n.$  Set $p=p(n,k)$ and $q=q(n,k).$ Then, by basic properties of the trace it follows that
\begin{align}
&\tau_{\frac{q}{p+sq}}\Big( \sigma_n(1_s\otimes x)y\Big)   \notag     \\ 
&= \tau_{\frac{q}{p+sq}}\Big( \sigma_n(1_s\otimes x)\sum_{i=1}^s \sigma_n(e_i\otimes 1_{\ga_n})y\sigma_n(e_i\otimes 1_{\ga_n})\Big) \notag \\
&= \frac{sq}{p+sq}\psi_s\otimes \tau_{\frac{p}{q}}\Big( \sigma_n^{-1}\Big(  \sigma_n(1_s\otimes x)\sum_{i=1}^s \sigma_n(e_i\otimes 1_{\ga_n})y\sigma_n(e_i\otimes 1_{\ga_n}) \Big)\Big) \notag \\
&= \frac{sq}{p+sq}\psi_s\otimes \tau_{\frac{p}{q}}\Big( (1_s\otimes x)V_n(y)\Big) \notag \\
&= \frac{sq}{p+sq}\tau_{\frac{p}{q}}\Big( x(\psi_s\otimes id_{\ga_n}(V_n(y)))\Big) \label{align:tracialrelate}
\end{align}
We now let $\widetilde{V}_s:\ga/\mathcal{J}_s\rightarrow \ga$ be the inductive limit of the maps 
$(\psi_s\otimes id_{\ga_n})\circ V_n$, which is well-defined by (\ref{eq:commutingsquare2}).  We also let 
$\widetilde{\sigma}_s:\ga\rightarrow \ga/\mathcal{J}_s$ be the inductive limit of the maps $\sigma_n\circ (1_s\otimes id_{\ga_n})$, which again are well-defined by (\ref{eq:commutingsquare2}).

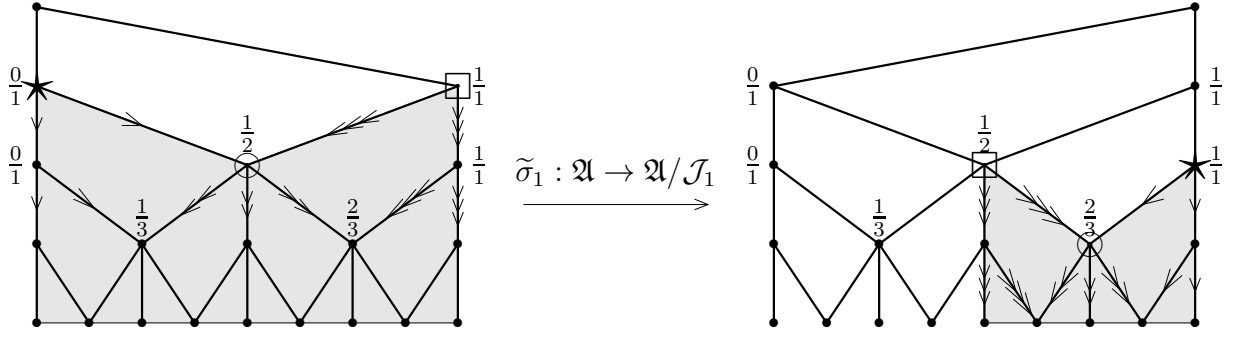
\begin{figure}[htb]
\begin{center}
\unitlength 0.35mm
\begin{picture}(400,120)(0,0)
\thinlines
\texture{c 0000}
\shade\path(-20,0)(-20,90)(60,60)(140,90)(140,0)(-20,0)
\shade\path(340,0)(420,0)(420,60)(380,30)(340,60)(340,0)

\linethickness{10mm} \path(230,47)(235,45)(230,43)
\path(-22,48)(-20,43)(-18,48) \path(12.5,75)(20,75)(14,80)
\path(-22,78)(-20,73)(-18,78) \path(408,47)(400,45)(405,52)
\path(410.5,20.5)(408,12)(415,18) \path(422,17)(420,12)(418,17)
\path(422,47)(420,42)(418,47) \path(-4,51)(2,43.5)(-6,47)

\path(40,47)(35,41)(42,44) \path(46,52)(41,46)(48,49)
\path(80,47)(85,41)(78,44) \path(74,52)(79,46)(72,49)
\path(58,45)(60,40)(62,45) \path(58,50)(60,45)(62,50)

\path(378,12)(380,7)(382,12) \path(378,17)(380,12)(382,17)
\path(367,15.5)(365,7.2)(371.5,13)
\path(371,21.4)(369,13.1)(375.5,18.9)
\path(389,21.4)(391,13.1)(384.5,18.9)
\path(393,15.5)(395,7.2)(388.5,13)

\path(97,76.2)(90,71.2)(98,72.2)
\path(103,78.5)(96,73.5)(104,74.5)
\path(109,80.8)(102,75.8)(110,76.8)

\path(119,47)(113,40)(121,43) \path(124,50.7)(118,43.7)(126,46.7)
\path(129,54.4)(123,47.4)(131,50.4)

\path(138,72)(140,67)(142,72) \path(138,77)(140,72)(142,77)
\path(138,82)(140,77)(142,82)

\path(138,42)(140,37)(142,42) \path(138,47)(140,42)(142,47)
\path(138,52)(140,47)(142,52)

\path(338,12)(340,7)(342,12) \path(338,17)(340,12)(342,17)
\path(338,22)(340,17)(342,22)

\path(338,42)(340,37)(342,42) \path(338,47)(340,42)(342,47)
\path(338,52)(340,47)(342,52)

\path(349.5,20.5)(352,12)(345,18)
\path(352.5,16)(355,7.5)(348,13.5) \path(355.5,11.5)(358,3)(351,9)

\path(352,47.2)(360,45.2)(355,52.2)\path(357,43.3)(365,41.3)(360,48.3)
\path(362,39.3)(370,37.3)(365,44.3)

\thicklines

\path(-20,0)(-20,120)(140,90)(140,0) \path(-20,90)(60,60)(140,90)
\path(60,0)(60,60) \path(20,0)(20,30) \path(100,0)(100,30)
\path(-20,60)(20,30)(60,60)(100,30)(140,60)

\path(260,0)(260,90)(420,120)(420,0) \path(260,90)(340,60)(420,90)
\path(340,0)(340,60) \path(300,0)(300,30) \path(380,0)(380,30)
\path(260,60)(300,30)(340,60)(380,30)(420,60)
\path(-20,30)(0,0)(20,30)(40,0)(60,30)(80,0)(100,30)(120,0)(140,30)
\path(260,30)(280,0)(300,30)(320,0)(340,30)(360,0)(380,30)(400,0)(420,30)

%\put(57,50){\makebox(0,0){$\frac{p^\prime}{q^\prime}$}}
%\put(57,0){\makebox(0,0){$\frac{p^\prime}{q^\prime}$}}

\put(-20,0){\makebox(0,0){\tiny $\bullet$}}
\put(0,0){\makebox(0,0){\tiny $\bullet$}}
\put(20,0){\makebox(0,0){\tiny $\bullet$}}
\put(40,0){\makebox(0,0){\tiny $\bullet$}}
\put(60,0){\makebox(0,0){\tiny $\bullet$}}
\put(80,0){\makebox(0,0){\tiny $\bullet$}}
\put(100,0){\makebox(0,0){\tiny $\bullet$}}
\put(120,0){\makebox(0,0){\tiny $\bullet$}}
\put(140,0){\makebox(0,0){\tiny $\bullet$}}

\put(-20,30){\makebox(0,0){\tiny $\bullet$}}
\put(20,30){\makebox(0,0){\tiny $\bullet$}}
\put(60,30){\makebox(0,0){\tiny $\bullet$}}
\put(100,30){\makebox(0,0){\tiny $\bullet$}}
\put(140,30){\makebox(0,0){\tiny $\bullet$}}
\put(260,30){\makebox(0,0){\tiny $\bullet$}}
\put(300,30){\makebox(0,0){\tiny $\bullet$}}
\put(340,30){\makebox(0,0){\tiny $\bullet$}}
\put(380,30){\makebox(0,0){\large $\odot$}}
\put(420,30){\makebox(0,0){\tiny $\bullet$}}

\put(-20,60){\makebox(0,0){\tiny $\bullet$}}
\put(60,60){\makebox(0,0){\large $\odot$}}
\put(140,60){\makebox(0,0){\tiny $\bullet$}}
\put(260,60){\makebox(0,0){\tiny $\bullet$}}
\put(340,60){\makebox(0,0){\large $\boxdot$}}
\put(420,60){\makebox(0,0){\huge $\star$}}
\put(-20,90){\makebox(0,0){\huge $\star$}}
\put(140,90){\makebox(0,0){\large $\boxdot$}}
\put(260,90){\makebox(0,0){\tiny $\bullet$}}
\put(420,90){\makebox(0,0){\tiny $\bullet$}}
\put(-20,120){\makebox(0,0){\tiny $\bullet$}}
\put(420,120){\makebox(0,0){\tiny $\bullet$}}

\put(260,0){\makebox(0,0){\tiny $\bullet$}}
\put(280,0){\makebox(0,0){\tiny $\bullet$}}
\put(300,0){\makebox(0,0){\tiny $\bullet$}}
\put(320,0){\makebox(0,0){\tiny $\bullet$}}
\put(340,0){\makebox(0,0){\tiny $\bullet$}}
\put(360,0){\makebox(0,0){\tiny $\bullet$}}
\put(380,0){\makebox(0,0){\tiny $\bullet$}}
\put(400,0){\makebox(0,0){\tiny $\bullet$}}
\put(420,0){\makebox(0,0){\tiny $\bullet$}}

\put(-28,90){\makebox(0,0){$\frac{0}{1}$}}
\put(-28,60){\makebox(0,0){$\frac{0}{1}$}}

\put(148,90){\makebox(0,0){$\frac{1}{1}$}}
\put(148,60){\makebox(0,0){$\frac{1}{1}$}}

\put(60,72){\makebox(0,0){$\frac{1}{2}$}}
\put(20,40){\makebox(0,0){$\frac{1}{3}$}}
\put(100,40){\makebox(0,0){$\frac{2}{3}$}}

\put(252,90){\makebox(0,0){$\frac{0}{1}$}}
\put(252,60){\makebox(0,0){$\frac{0}{1}$}}

\put(428,90){\makebox(0,0){$\frac{1}{1}$}}
\put(428,60){\makebox(0,0){$\frac{1}{1}$}}

\put(340,72){\makebox(0,0){$\frac{1}{2}$}}
\put(300,40){\makebox(0,0){$\frac{1}{3}$}}
\put(380,40){\makebox(0,0){$\frac{2}{3}$}}

\put(200,57){\makebox(0,0){$\widetilde{\sigma}_1:\ga \rightarrow
\ga/\mathcal{J}_1$}}

\thinlines \path(165,45)(235,45)

\end{picture}
\end{center}
\caption{The map $\widetilde{\sigma}_1$} \label{Figure3}
\end{figure}

Figure \ref{Figure3} (graciously provided by F. Boca) displays the mapping $\widetilde{\sigma}_1$ in terms of the Bratelli diagrams of $\ga$ and $\ga/\mathcal{J}_1.$

Set $\pi=\oplus_{\frac{1}{s+1}<\frac{p}{q}<\frac{1}{s}}\pi_{\frac{p}{q}}$ and identify $\ga/\mathcal{J}_s$ with $\pi(\ga)$ by (\ref{eq:anotherdescriptionJ_s}).
By the Choi-Effros lifting theorem in \cite{Choi76}, there is a UCP lifting $\phi:\ga/\mathcal{J}_s\rightarrow \ga$ of $\pi.$  Then let
\begin{equation*}
 G_s:=\widetilde{V}_s\circ \pi :\ga\rightarrow\ga\quad \textrm{ and }H_s:=\phi\circ \widetilde{\sigma}_s:\ga\rightarrow \ga.
\end{equation*}
By (\ref{eq:paritalinversesonfd}), it follows that
\begin{equation*}
 G_sH_s=id_\ga.
\end{equation*}
It is also routine to verify, using the definitions of $\sigma_n$ and $V_n$, that
\begin{equation}
 G_s(xH_s(y))=G_s(x)y\quad \textrm{ for every }x,y\in \ga.  \label{eq:shiftingGandH}
\end{equation}

By (\ref{align:tracialrelate}), we have the following relationship for every $x,y\in\ga,$ and $\frac{p}{q}\in\mathbb{Q}\cap[0,1]:$
\begin{align*}
&=\tau_{g_s(p/q)}^\ga(H_s(x)y)=\tau_{\frac{q}{p+sq}}^\ga(H_s(x)y)=\tau_{\frac{q}{p+sq}}\Big( \pi_{\frac{q}{p+sq}}(H_s(x)y)\Big)\\
&=\frac{sq}{p+sq}\tau_{\frac{p}{q}}\Big( \pi_{\frac{p}{q}}(xG_s(y))\Big)=\frac{sq}{p+sq}\tau_{\frac{p}{q}}^\ga(xG_s(y))=sg_s(p/q)\tau^\ga_{\frac{p}{q}}(xG_s(y))
\end{align*}

Therefore, by Theorem \ref{thm:statetrace}, for any $\theta\in [0,1]$ we have
\begin{equation}
 \tau_{g_s(\theta)}^\ga(H_s(x)y)= sg_s(\theta)\tau_{\theta}^\ga(xG_s(y))\label{eq:maintracerelationshipall}
\end{equation}

Therefore, by Corollary \ref{cor:traceideal}, it follows that for any $\theta\in[0,1]$, we have
\begin{equation}
 G_s(I_{g_s(\theta)})=I_{\theta}.  \label{eq:idealcorrespondence}
\end{equation}

Moreover, by the description of $Z(\ga)$ given in Theorem \ref{thm:conditionalexpectationdefn} and (\ref{eq:maintracerelationshipall}) it is clear that
\begin{equation}
G_s(f)=f\circ g_s\quad \textrm{ for every }f\in C[0,1]. \label{eq:Gsandcomposition}
\end{equation}
For each $s\geq1,$ define $f_s\in Z(\ga)\cong C[0,1]$ as 
\begin{equation}
 f_s(\theta)=\frac{\theta+1}{(\theta+s)(\theta+s+1)}  \label{eq:defnofthef_s}
\end{equation}
Let us now define $\widetilde{\mathbb{G}}:\ga\rightarrow \ga$ as 
\begin{equation}
 \widetilde{\mathbb{G}}(x)=\sum_{s=1}^\infty G_s(x)f_s.
\end{equation}
\section{Proof of Theorem \ref{thm:maintheorem}} 
In this section we will prove the 5 assertions from Theorem \ref{thm:maintheorem}. First note that Theorem \ref{thm:maintheorem}(1) follows from (\ref{eq:V_G*definition}) and (\ref{eq:Gsandcomposition}), and (2) follows from (\ref{eq:idealcorrespondence}).

First define $\phi_0:=\tau_0^\ga.$  Then, let $\theta\in [0,1]$ with $\frac{1}{s+1}<\theta\leq \frac{1}{s}$ for some $s\geq1.$  Then define 
\begin{equation*}
 \phi_\theta(x)=\tau_\theta^\ga(H_s(1))^{-1}\tau_\theta^\ga(H_s(1)x)=\frac{1}{s\theta}\tau_\theta^\ga(H_s(1)x)  \quad \textrm{ for every }\quad x\in\ga.
\end{equation*}
Recall that Gauss measure $\mu$ on $[0,1]$ is defined as the probability measure $d\mu=\frac{d\theta}{\ln 2(\theta+1)}$, where $d\theta$ is Lebesgue measure. Let $\phi\in \mathcal{S}(\ga)$ be the direct integral of the states $\phi_\theta$ over  $\mu$, i.e.
\begin{equation*}
 \phi(x)=\int_0^1 \phi_\theta(x)d\mu(\theta).
\end{equation*}
Let $\tau\in\mathcal{T}(\ga)$ be the unique tracial extension of $\mu$ provided by Theorem \ref{thm:statetrace}.  By uniqueness we have
\begin{equation*}
 \tau=\int_0^1 \tau^\ga_\theta d\mu(\theta).
\end{equation*}

Notice that for every $f\in C[0,1]$ and $x\in\ga$, we have 
\begin{equation}
\phi_\theta(fx)=f(\theta)\phi_\theta(x)\quad \textrm{and}\quad\tau_\theta^\ga(fx)=f(\theta)\phi_\theta(x) \label{eq:Gaussmeasurerestriction}
\end{equation}
It also follows from (\ref{eq:Gaussmeasurerestriction}) that $\phi$ restricted to $C[0,1]$ is Gauss measure $\mu.$

For any state $\psi\in \mathcal{S}(\ga)$, let $(L^2(\ga,\psi), \pi_\psi)$ denote the GNS representation of $\psi$ and $\la\cdot,\cdot\ra_{\psi}$ the inner product on $L^2(\ga,\psi).$  For $x\in \ga$, we will denote by $x_\psi$ the image of $x$ in $L^2(\ga,\psi)$ and denote by $\ga_\psi$ the dense subspace of $L^2(\ga,\psi)$ consisting of  the $x_\psi.$

By the definitions of $\phi$ and $\tau$, we can decompose 
\begin{equation*}
 L^2(\ga,\phi)=\int_0^1 L^2(\ga,\phi_\theta)d\mu(\theta)\quad \textrm{and}\quad L^2(\ga,\tau)=\int_0^1 L^2(\ga,\tau_\theta^\ga)d\mu(\theta). 
\end{equation*}

Furthermore, by (\ref{eq:Gaussmeasurerestriction}) we have
\begin{equation*}
 L^2(\mu)\subset L^2(\ga,\phi)\quad \textrm{and}\quad L^2(\mu)\subset L^2(\ga,\tau)
\end{equation*}
as
\begin{equation}
f_\phi=\int_0^1 f(\theta)1_{\phi_\theta}d\mu(\theta) \quad \textrm{and}\quad f_\tau=\int_0^1 f(\theta)1_{\tau^\ga_\theta}d\mu(\theta). \label{eq:commuteinsidenoncommute}
\end{equation}

We now define an isometry $\widetilde{V}_G:L^2(\ga,\tau)\rightarrow L^2(\ga,\phi)$ that satisfies (3)-(5) in Theorem \ref{thm:maintheorem}.
As short hand notation, for each vector $\eta\in L^2(\ga,\phi)$ and Borel set $E\subset [0,1]$ we will write 
\begin{equation*}
 \eta 1_E:=\int_E \eta(\theta)d\mu(\theta)\in  \int_0^1 L^2(\ga,\phi_\theta)d\mu(\theta)
\end{equation*}

For each $s\geq1$, define operators on $\ga_\phi$ and $\ga_\tau$ respectively as
\begin{equation}
 \widetilde{H}_s(x_\tau)=H_s(x)_\phi 1_{[\frac{1}{s},\frac{1}{s+1}]}\quad \textrm{and }\widetilde{G}_s(x_\phi)=(G_s(x)f_s)_\tau.
\end{equation}
Clearly these maps are contractive, so they extend to operators on $L^2(\ga,\phi)$ and $L^2(\ga,\tau)$ respectively.  Now define 
\begin{equation*}
 \widetilde{V}_G=WOT-\sum_{s=1}^\infty \widetilde{H}_s.
\end{equation*}
We now show that $\widetilde{V}_G$ is an isometry.  Let us first recall $f_s$ from (\ref{eq:defnofthef_s}) and note that 
\begin{equation*}
 \sum_{s=1}^\infty f_s(\theta)=1 \quad \textrm{ for every }\theta\in[0,1].
\end{equation*}
We will implicitly use this fact throughout the rest of the proof of Theorem \ref{thm:maintheorem}. We have,

\begin{align}
 \la \widetilde{V}_G(x_\tau),\widetilde{V}_G(x_\tau)\ra_{\phi}&= \sum_{s=1}^\infty\int_{\frac{1}{s+1}}^{\frac{1}{s}}\phi_\theta(H_s(x)^*H_s(x))d\mu(\theta) \notag \\
&=\sum_{s=1}^\infty\int_{\frac{1}{s+1}}^{\frac{1}{s}} \frac{1}{s\theta}\tau^\ga_\theta(H_s(x^*x))d\mu(\theta)\notag \\
&=\sum_{s=1}^\infty\int_{\frac{1}{s+1}}^{\frac{1}{s}} \tau^\ga_{\frac{1}{\theta}-s}(x^*x)d\mu(\theta)\quad (\textup{by (\ref{eq:maintracerelationshipall})}) \notag \\
&=\frac{1}{\ln 2}\sum_{s=1}^\infty\int_0^1 \tau_u^\ga(x^*x)\frac{u+1}{(u+s)(u+s+1)}\frac{du}{u+1}\quad (\textup{with }u=\frac{1}{\theta}-s) \label{align:V_Gisometry} \\
&=\int_0^1\tau_u^\ga(x^*x)\Big(\sum_{s=1}^\infty f_s(u)  \Big)d\mu(u) \notag \\
&=\la x_\tau,x_\tau\ra_{\tau}. \notag
\end{align}

We now calculate $\widetilde{V}_G^*.$  Let $x,y\in\ga,$ then 
\begin{align*}
\la \widetilde{V}_G(x_\tau),y_\phi\ra_\phi&=\sum_{s=1}^\infty \int_{\frac{1}{s+1}}^{\frac{1}{s}}\phi_\theta(y^*H_s(x))d\mu(\theta)\\
&=\sum_{s=1}^\infty \int_{\frac{1}{s+1}}^{\frac{1}{s}}\frac{1}{s\theta}\tau_\theta^\ga(y^*H_s(x))d\mu(\theta)\\
&=\sum_{s=1}^\infty \int_{\frac{1}{s+1}}^{\frac{1}{s}}\tau_{\frac{1}{\theta}-s}^\ga(G_s(y)^*x)d\mu(\theta)\quad (\textup{by (\ref{eq:maintracerelationshipall})})\\
&=\sum_{s=1}^\infty \int_0^1\tau_\theta^\ga(G_s(y)^*x)f_s(\theta)d\mu(\theta)\quad (\textup{Reasoning as in (\ref{align:V_Gisometry})})\\
&=\sum_{s=1}^\infty \int_0^1\tau_\theta^\ga(G_s(y)^*f_s x)d\mu(\theta)\quad (\textup{By (\ref{eq:Gaussmeasurerestriction})})\\
&=\la x_\tau, \sum \widetilde{G}_s(y_\phi)\ra_\tau.\\
&=\la x_\tau,\widetilde{\mathbb{G}}(y)_\tau\ra_\tau.
\end{align*}
We now show (3).  Let $f=\int_0^1 f(\theta)1_{\tau_\theta}d\mu(\theta)\in L^2(\mu)\subset L^2(\ga,\tau).$ Then

\begin{align}
\widetilde{V}_G(f)&=\sum_{s=1}^\infty\int_{\frac{1}{s+1}}^{\frac{1}{s}}f(\frac{1}{\theta}-s)(H_s(1))_{\phi_\theta}d\mu(\theta) \notag \\
&=\sum_{s=1}^\infty\int_{\frac{1}{s+1}}^{\frac{1}{s}}f(\frac{1}{\theta}-s)1_{\phi_\theta}d\mu(\theta)\notag \\
&=f\circ G \in L^2(\ga,\phi).\label{align:V_goncommutative}
\end{align}
Similarly, one shows that $\widetilde{V}_G^*|_{L_2(\mu)}=V_G^*$.  This proves (3).

We now show (4).  It follows from the definition of $G_s$ that for every $s\geq1$ we have 
\begin{equation*}
 \int_{[\frac{1}{s+1},\frac{1}{s}]^c}L^2(\ga,\phi_\theta)d\mu(\theta)\subset \textrm{ker}(\widetilde{G}_s)
\end{equation*}
From this and (\ref{eq:shiftingGandH}) it follows that for every $x,y\in\ga$ we have
\begin{align*}
 \widetilde{V}_G^*\pi_\phi(x)\widetilde{V}_G(y_\tau)&=\widetilde{V}_G^*\Big(\sum_{s=1}^\infty (xH_s(y))_\phi 1_{[\frac{1}{s+1},\frac{1}{s}]}\Big)\\
&=\sum_{s=1}^\infty(G_s(xH_s(y))f_s)_\tau\\
&=\sum_{s=1}^\infty(G_s(x)yf_s)_\tau\quad (\textup{By (\ref{eq:shiftingGandH})})\\
&=\pi_\tau(\widetilde{\mathbb{G}}(x))y_\tau.
\end{align*}

By (\ref{align:V_goncommutative}) we have $\widetilde{V}_G(1_\tau)=1_\phi$, from which it follows that
\begin{equation*}
 \phi(x)=\la x_\phi, 1_\phi\ra =\la x_\phi,\widetilde{V}_G(1_\tau)\ra =\la \widetilde{\mathbb{G}}(x)_\tau,1_\tau\ra=\tau(\widetilde{\mathbb{G}}(x)).
\end{equation*}
This proves (5) and finishes the proof of Theorem \ref{thm:maintheorem}.
\section*{Acknowledgment}
Many thanks are due to Florin Boca for suggesting this problem to me, several helpful conversations, and for providing code for the figures in this article.

%\section{Sketch of map into multiplier algebra}
%Let $x\in\ga.$  Since the set $\{ \tau_\theta^\ga: \theta\in[0,1]\}$ equipped with the $\sigma(\ga^*, \ga)$ topology is homeomorphic to $[0,1]$ it follows that the function defined by
%\begin{equation*}
 %f_x(\frac{p}{q})=\tau^\ga_{\frac{q}{p+sq}}(x)\quad \textrm{ for }\frac{p}{q}\in\mathbb{Q}\cap[0,1]
%\end{equation*}
%extends to a continuous function on $[0,1]$, which we still denote by $f_x.$  Consider $H_s$ as a map into $\ga/\mathcal{J}_s\subset (\ga/\mathcal{J}_s)^{**}\subset\ga^{**}.$ Let $\pi:\ga\rightarrow \ga/\mathcal{J}_s$ be the quotient map. Consider the map $\widetilde{H}_s:\ga\rightarrow \ga^{**}$ defined by
%\begin{equation*}
 %\widetilde{H}_s=H_s(x)+ \pi_s\Big((1-H_s(1))f_x\Big). \quad \textrm{ for every }x\in \ga
%\end{equation*}
%Let 
%\begin{equation*}
 %\mathbb{G}(x)=WOT-\sum_{s=1}^\infty H_s(x)\in\ga^{**}. \quad\textrm{ for }x\in\ga.
%\end{equation*}
%By the definition of $\widetilde{H}_s$  it is clear that
%\begin{equation*}
 %\mathbb{G}(f)=f\circ G\in C^b(C_0(0,1)) \quad \textrm{ for }f\in C_0(0,1)\subset \ga.
%\end{equation*}
%It is at least plausible and almost certainly true that $\mathbb{G}$ maps $\mathcal{J}_{\{0,1\}}$ into $\mathcal{M}(\mathcal{J}_{\{0,1\}}).$

\bibliographystyle{plain}
\bibliography{mybib}

\begin{thebibliography}{1}

\bibitem{Boca08}
Florin~P. Boca.
\newblock An {AF} algebra associated with the {F}arey tessellation.
\newblock {\em Canad. J. Math.}, 60(5):975--1000, 2008.

\bibitem{Choi76}
Man~Duen Choi and Edward~G. Effros.
\newblock The completely positive lifting problem for {$C\sp*$}-algebras.
\newblock {\em Ann. of Math. (2)}, 104(3):585--609, 1976.

\bibitem{Effros80}
Edward~G. Effros and Chao~Liang Shen.
\newblock Approximately finite {$C\sp{\ast} $}-algebras and continued
  fractions.
\newblock {\em Indiana Univ. Math. J.}, 29(2):191--204, 1980.

\bibitem{Iosifescu02}
Marius Iosifescu and Cor Kraaikamp.
\newblock {\em Metrical theory of continued fractions}, volume 547 of {\em
  Mathematics and its Applications}.
\newblock Kluwer Academic Publishers, Dordrecht, 2002.

\bibitem{Mundici88}
Daniele Mundici.
\newblock Farey stellar subdivisions, ultrasimplicial groups, and {$K\sb 0$} of
  {AF} {$C\sp *$}-algebras.
\newblock {\em Adv. in Math.}, 68(1):23--39, 1988.

\bibitem{Mundici08}
Daniele Mundici.
\newblock Revisiting the {F}arey {AF} algebra.
\newblock preprint, 2008.

\bibitem{Pedersen79}
Gert~K. Pedersen.
\newblock {\em ${C}^*$-algebras and their automorphism groups.}
\newblock Number~14 in London Mathematical Society Monographs. Academic Press
  Inc., London-New York, 1979.

\end{thebibliography}

\end{document}